\documentclass[11pt]{article}

\usepackage{amsmath}
\usepackage{amsthm}
\usepackage{amssymb}
\usepackage{latexsym}
\usepackage{pstricks}
\usepackage{graphicx, pst-plot, pst-node, pst-text, pst-tree}
\usepackage{titlefoot}
\usepackage{enumerate}
\usepackage{titlesec}
\usepackage[small,it]{caption}
\usepackage{color}
\definecolor{lightgray}{rgb}{0.8, 0.8, 0.8}
\definecolor{darkgray}{rgb}{0.7, 0.7, 0.7}
\definecolor{darkblue}{rgb}{0, 0, .4}

\usepackage[bookmarks,breaklinks]{hyperref}
\hypersetup{
        colorlinks=true,
        linkcolor=darkblue,
        anchorcolor=darkblue,
        citecolor=darkblue,
        pagecolor=darkblue,
        urlcolor=darkblue,
        pdfpagemode=UseThumbs,
        pdftitle={Reconstructing compositions},
        pdfsubject={Combinatorics},
        pdfauthor={Vincent Vatter},
}

\newtheorem{theorem}{Theorem}

\newtheorem{lemma}[theorem]{Lemma}

\newtheoremstyle{example}{\topsep}{\topsep}%
     {}%         Body font
     {}%         Indent amount (empty = no indent, \parindent = para indent)
     {\bfseries}% Thm head font
     {.}%        Punctuation after thm head
     {.5em}%     Space after thm head (\newline = linebreak)
     {\thmname{#1}\thmnumber{ #2}}%         Thm head spec
\theoremstyle{example}
\newtheorem{example}[theorem]{Example}

\newcounter{todocounter}

\newcommand{\minisec}[1]{\bigskip\noindent{\bf #1.}}

\long\def\symbolfootnote[#1]#2{\begingroup%
\def\thefootnote{\fnsymbol{footnote}}\footnote[#1]{#2}\endgroup}

\newcommand{\ex}{\operatorname{ex}}
\newcommand{\getwo}{\mathord{\ge}2}

\makeatletter
\newfont{\footsc}{cmcsc10 at 8truept}
\newfont{\footbf}{cmbx10 at 8truept}
\newfont{\footrm}{cmr10 at 10truept}
\renewenvironment{abstract}%
                {
                  \begin{list}{}%
                     {\setlength{\rightmargin}{1in}%
                      \setlength{\leftmargin}{1in}}%
                   \item[]\ignorespaces\begin{small}}%
                 {\end{small}\unskip\end{list}}

\datefoot{\today}
\amssubj{68R15, 06A07, 05A17, 05C60}

\newpagestyle{main}[\small]{
        \headrule
        \sethead[\usepage][][]
        {\sc Reconstructing Compositions}{}{\usepage}}

\title{\sc{Reconstructing Compositions}}
\author{\sc{Vincent Vatter\thanks{Supported by EPSRC grant GR/S53503/01.}}\\
\small School of Mathematics and Statistics\\[-3pt]
\small University of St Andrews\\[-3pt]
\small St Andrews, Fife, Scotland\\[-3pt]
\small \texttt{vince@mcs.st-and.ac.uk}\\[-3pt]
\small \href{http://turnbull.mcs.st-and.ac.uk}{\texttt{http://turnbull.mcs.st-and.ac.uk/\~{}vince}}\\[-10pt]}

\titleformat{\section}
        {\large\sc}
        {\thesection.}{1em}{}   

\date{}

\begin{document}
\maketitle

\pagestyle{main}

\begin{abstract}
We consider the problem of reconstructing compositions of an integer from their subcompositions, which was raised by Raykova (albeit disguised as a question about layered permutations).  We show that every composition $w$ of $n\ge 3k+1$ can be reconstructed from its set of $k$-deletions, i.e., the set of all compositions of $n-k$ contained in $w$.  As there are compositions of $3k$ with the same set of $k$-deletions, this result is best possible.
\end{abstract}

\minisec{Introduction}
The Reconstruction Conjecture states that given the multiset of isomorphism types of $1$-vertex deletions (briefly, {\it $1$-deletions\/}) of a graph $G$ --- the {\it deck\/} of $G$ ---  on three or more vertices, it is possible to determine $G$ up to isomorphism.  The stronger set version of the conjecture due to Harary~\cite{harary:on-the-reconstr:} only allows access to the {\it set\/} of $1$-deletions and requires $G$ to have four or more vertices.  These conjectures can be made even more difficult by considering $k$-deletions instead of $1$-deletions, for which we refer to Manvel~\cite{manvel:some-basic-obse:}.

Such reconstruction questions extend naturally to other combinatorial contexts.  For example, Sch\"utzenberger and Simon (see Lothaire~\cite[Theorem 6.2.16]{lothaire:combinatorics-o:}) proved that every word of length $n\ge 2k+1$ can be reconstructed from its set of $k$-deletions (i.e., subwords of length $n-k$).  This bound is tight because the words $(ab)^k$ (the word with $ab$ repeated $k$ times) and $(ba)^k$ have the same set of $k$-deletions: all words of length $k$ over the set $\{a,b\}$.  Answering a question of Cameron~\cite{cameron:stories-from-th:}, Pretzel and Siemons~\cite{pretzel:on-the-reconstr:} considered the partition context, where they proved that every partition of $n\ge 2(k+3)(k+1)$ can be reconstructed from its set of $k$-deletions.  (This bound is not known to be tight.)

Motivated by a question of Raykova~\cite{raykova:permutation-rec:} (described at the end of the paper), we consider the problem of set reconstruction for compositions (ordered partitions), establishing the following result.

\begin{theorem}\label{thm-comp-reconstruct}
All compositions of $n\ge 3k+1$ can be reconstructed from their sets of $k$-deletions.
\end{theorem}

Our proof of Theorem~\ref{thm-comp-reconstruct} illustrates an algorithm to perform the reconstruction.  Perhaps more convincing than the proof is the Maple implementation of this algorithm, available from the author's homepage.

\minisec{Notation}
We view a composition as a word $w$ whose letters are positive integers, i.e., a word in $\mathbb{P}^*$.  We denote the length of $w$ by $|w|$ and the sum of the entries of $w$ by $\|w\|$, and say that $w$ is a composition of $\|w\|$.  A $1$-deletion of $w$ is a composition that can be obtained either by lowering a $\getwo$ entry of $w$ by $1$ or by removing an entry of $w$ that is equal to $1$.  A $2$-deletion is then a $1$-deletion of a $1$-deletion, and so on.

This notion naturally defines a partial order%
\symbolfootnote[2]{This partial order was first considered by Bergeron, Bousquet-M\'elou, and Dulucq~\cite{bergeron:standard-paths-:}, and has since been studied by Snellman~\cite{snellman:saturated-chain:,snellman:standard-paths-:}, Sagan and Vatter~\cite{sagan:the-mobius-func:}, and Bj\"orner and Sagan~\cite{bjorner:rationality-of-:}.}
on compositions: $u\le w$ if $w$ contains a subword $w(i_1)w(i_2)\cdots w(i_\ell)$ of length $\ell=|u|$ such that $u(j)\le w(i_j)$ for all $1\le j\le \ell$.  (We refer to the indices $i_1<\cdots<i_\ell$ as an {\it embedding\/} of $u$.)  For example, $1211\le 21312$ because of the subword $2312$.  If $u\le w$ then $u$ is a $(\|w\|-\|u\|)$-deletion of $w$.  Returning to the previous example, $\|21312\|=9$ and $\|1211\|=5$, so $1211$ is a $4$-deletion of $21312$.

\minisec{A lower bound}
In the context of words, the fact that the sets of $k$-deletions of $(ab)^k$ and $(ba)^k$ are both equal to the set of all words of length $k$ over $\{a,b\}$ provides a lower bound on $k$-reconstructibility.  Here we can use a very similar example: the sets of $k$-deletions of $(12)^k$ and $(21)^k$ are both equal to the set of all compositions of $2k$ in which no entry is greater than $2$.  This implies that Theorem~\ref{thm-comp-reconstruct} is best possible.

\minisec{The proof}
Our reconstruction algorithm/proof of Theorem~\ref{thm-comp-reconstruct} employs several composition statistics.  One is the {\it exceedance number\/}, defined by $\ex(w)=\|w\|-|w|=\sum (w_i-1)$ where the sum is over all entries $w(i)$.  Another important composition statistic is the number of $1$'s in $w$, which can be approximated using its set of $k$-deletions:

\begin{lemma}\label{lem-determine-ones}
The composition $w$ of $n\ge 3k+1$ has at least $k$ $1$'s if and only if either
\begin{enumerate}
\item[(1)] $1^{n-k}$ is a $k$-deletion of $w$, or
\item[(2)] the longest $k$-deletion of $w$ is $k$ letters longer than the shortest $k$-deletion of $w$.
\end{enumerate}
Moreover, $w$ has precisely $k$ $1$'s if and only if one of the above conditions holds and $w$ has a $k$-deletion without $1$'s.
\end{lemma}
\begin{proof}
It is easy to see that if either (1) or (2) occurs then $w$ has at least $k$ $1$'s.  Suppose then that $w$ has at least $k$ $1$'s.  If $\ex(w)\le k$ then $1^{n-k}$ is a $k$-deletion of $w$, satisfying (1).  On the other hand, if $\ex(w)>k$ then some $k$-deletion of $w$ has length $|w|$, while the fact that $w$ contains at least $k$ $1$'s guarantees that some $k$-deletion of $w$ has length $|w|-k$, satisfying (2).  The second claim in the lemma is then readily verified.
\end{proof}

Given a set of $k$-deletions of a composition, the first step in our algorithm is to apply Lemma~\ref{lem-determine-ones} to decide if the composition has fewer than $k$, precisely $k$, or more than $k$ $1$'s.  The three cases are handled separately.  The first two are relatively straightforward, while the last is more delicate.

\begin{lemma}\label{lem-few-ones}
If $w$ is a composition of $n\ge 3k+1$ with fewer than $k$ $1$'s, then $w$ can be reconstructed from its set of $k$-deletions.
\end{lemma}
\begin{proof}
Given the set of $k$-deletions of a composition $w$ satisfying these hypotheses, our algorithm can apply the result of Lemma~\ref{lem-determine-ones} to determine that $w$ has fewer than $k$ $1$'s.  It then follows that
$$
\ex(w)\ge\frac{\|w\|-(\mbox{\# of $1$'s in $w$})}{2}\ge\frac{2k+2}{2}=k+1.
$$
From this we see that $w$ has the same length, say $m$, as its longest $k$-deletions, and then $\ex(w)$ can be easily determined: it is $k$ plus the exceedance number of one of the longest $k$-deletions.

Set $t=\ex(w)-k$ and define the composition $a=a(1)\cdots a(m)$ by
$$
a(i)=\max\{s : \mbox{$\underbrace{1\cdots1}_{i-1}s\underbrace{1\cdots1}_{m-i}$ is, or is contained in, a $k$-deletion of $w$}\}.
$$
It follows that $a$ satisfies
\begin{equation}\label{eqn-u-v-w}
a(i)=\min\{w(i), t+1\}.
\end{equation}
There are now two cases in which we are done:
\begin{itemize}
\item If $\|a\|=n$ then $w$ must be equal to $a$.  By \eqref{eqn-u-v-w}, this will occur if $w$ contains no entries greater than $t+1$.
\item If at most one entry of $a$ satisfies $a(i)=t+1$ --- which by \eqref{eqn-u-v-w} will occur if $w$ contains at most one entry $w(i)\ge t+1$ --- then \eqref{eqn-u-v-w} forces $w(j)=a(j)$ for all $j\neq i$ and then $w(i)$ can be calculated from the fact that $\|w\|=n$.
\end{itemize}
Suppose, for the sake of contradiction, that neither of these conditions hold.  Thus $w$ must contain an entry $w(i)>t+1$ and another entry $w(j)\ge t+1$.  We then have
$$
k+t=\ex(w)\ge t + (t+1) + (\mbox{\# of $\getwo$ entries in $w$, not including $w(i),w(j)$}),
$$
so
\begin{equation}\label{eqn-few-ones-2}
k\ge t+1+(\mbox{\# of $\getwo$ entries in $w$, not including $w(i),w(j)$}),
\end{equation}
while
$$
|w|=2+(\mbox{$\#1$s in $w$})+(\mbox{\# of $\getwo$ entries in $w$, not including $w(i),w(j)$}),
$$
so because $w$ contains fewer than $k$ $1$'s,
\begin{equation}\label{eqn-few-ones-3}
(\mbox{\# of $\getwo$ entries in $w$, not including $w(i),w(j)$})\ge |w|-k-1.
\end{equation}
Combining (\ref{eqn-few-ones-2}) and (\ref{eqn-few-ones-3}) shows that $|w|\le 2k-t$, but then $\ex(w)\ge (3k+1)-(2k-t)=k+t+1$, contradicting the definition of $t$ and completing the proof.
\end{proof}

\begin{example}
Suppose the reconstruction algorithm is given the set of $3$-deletions
$$
\{52,
322, 412, 421, 511,
2122, 3112, 3121, 4111\}
$$
of an unknown composition $w$ of $n=10$.  The algorithm first checks the hypotheses of Lemma~\ref{lem-determine-ones}.  The first condition does not hold because the set of $3$-deletions does not contain $1^{10-3}=1111111$, while the second condition fails because the longest $3$-deletion is only $2$ letters longer than the shortest.  Therefore $w$ has fewer than $k=3$ $1$'s.  Now the algorithm follows the proof of Lemma~\ref{lem-few-ones}.  First we compute $\ex(w)$ from one of the longest $3$-deletions:
$$
\ex(w)=\ex(3121)+3=6,
$$
so $t=3$.  Then we compute $a$:
$$
\begin{array}{l}
\mbox{$a(1)=4$ because $4111$ is contained in a $3$-deletion but $5111$ is not,}\\
\mbox{$a(2)=1$ because $1111$ is contained in a $3$-deletion but $1211$ is not,}\\
\mbox{$a(3)=2$ because $1121$ is contained in a $3$-deletion but $1131$ is not,}\\
\mbox{$a(4)=2$ because $1112$ is contained in a $3$-deletion but $1113$ is not.}
\end{array}
$$
Thus $w\ge 4122$.  Since $\|4122\|=9<10=\|w\|$, we are not done reconstructing $w$ and need to account for one more exceedance.  However, since $a(1)$ is the only entry of $a$ equal to $t+1=4$, $w(1)$ is the only entry of $w$ that can be greater than the corresponding entry of $a$, so we get $w=5122$.
\end{example}

\begin{lemma}\label{lem-k-ones}
If $w$ is a composition of $n\ge 3k+1$ with precisely $k$ $1$'s, then $w$ can be reconstructed from its set of $k$-deletions.
\end{lemma}
\begin{proof}
Given the set of $k$-deletions of a composition $w$ satisfying these hypotheses, our algorithm can apply the result of Lemma~\ref{lem-determine-ones} to determine that it has exactly $k$ $1$'s.  With this established, the length of $w$ can be computed as $k$ plus the length of the shortest $k$-deletion of $w$.

There is a $k$-deletion of $w$ without $1$'s, and this composition gives the $\getwo$ entries of $w$ in their correct order.  Thus it suffices to determine where they lie in $w$.  To this end define the composition $a_i$ by
$$
a_i=\underbrace{1\cdots1}_{i-1}2\underbrace{1\cdots1}_{m-i}.
$$
As $a_i$ is contained in a $k$-deletion of $w$ if and only if $w(i)\ge 2$, the $\getwo$ entries of $w$ can be discerned, completing the proof.
\end{proof}

\begin{example}
Suppose the reconstruction algorithm is given the set of $3$-deletions
$$
\begin{array}{l}
\{322,
2212, 2221, 3112, 3121, 3211,
12121, 12211, 21121,\\
\ \ 21211, 22111, 31111,
111211, 121111, 211111\}.
\end{array}
$$
of an unknown composition $w$ of $n=10$.  Since the longest $3$-deletions in this set are $3$ letters longer than the shortest $3$-deletion, $w$ has at least $k=3$ $1$'s by Lemma~\ref{lem-determine-ones}.  As the set also contains a $3$-deletion without $1$'s, the same lemma shows that $w$ has precisely $3$ $1$'s, and thus the algorithm follows the proof of Lemma~\ref{lem-k-ones}.  The $3$-deletion without $1$'s --- $322$ --- gives the $\getwo$ entries of $w$ in their correct order.  Now we form the $a_i$'s to see where these $\getwo$ entries lie:
$$
\begin{array}{l}
\mbox{$a_1=211111$ is contained in a $3$-deletion so $w(1)\ge 2$,}\\
\mbox{$a_2=121111$ is contained in a $3$-deletion so $w(2)\ge 2$,}\\
\mbox{$a_3=112111$ is not contained in a $3$-deletion so $w(3)=1$,}\\
\mbox{$a_4=111211$ is contained in a $3$-deletion so $w(4)\ge 2$,}\\
\mbox{$a_5=111121$ is not contained in a $3$-deletion so $w(5)=1$,}\\
\mbox{$a_6=111112$ is not contained in a $3$-deletion so $w(6)=1$.}
\end{array}
$$
Therefore we get $w=321211$.
\end{example}

This leaves us to consider the case of compositions with many $1$'s.  In this case we also need the {\it second exceedance number\/}, defined by $\ex_2(w)=\sum (w(i)-2)$ where the sum is over all entries $w(i)\ge 2$.

\begin{lemma}\label{lem-many-ones}
If $w$ is a composition of $n\ge 3k+1$ with more than $k$ $1$'s, then $w$ can be reconstructed from its set of $k$-deletions.
\end{lemma}
\begin{proof}
Given the set of $k$-deletions of such a composition $w$, our algorithm can apply the result of Lemma~\ref{lem-determine-ones} to conclude that it has more than $k$ $1$'s.  Therefore the $k$-deletions with the fewest $1$'s contain all $\getwo$ entries of $w$ in the order in which they occur in $w$; let $v=v(1)\cdots v(\ell)$ denote the composition formed by these entries, so
$$
w=\underbrace{1\cdots1}_{z(1)}v(1)\underbrace{1\cdots1}_{z(2)}v(2)\cdots v(\ell-1)\underbrace{1\cdots1}_{z(\ell)}v(\ell)\underbrace{1\cdots1}_{z(\ell+1)}
$$
for some word $z\in\mathbb{N}^{\ell+1}$ (we take $\mathbb{N}$ to denote the nonnegative integers).  Our goal is thus to determine $z$.  We use similar techniques as in the proof of Lemma~\ref{lem-few-ones}, although here we must perform two steps.

The first of these steps is to find the $0$'s in $z$.  For $1\le i\le\ell+1$ let
$$
a_i=\underbrace{2\cdots2}_{i-1}1\underbrace{2\cdots2}_{\ell+1-i}.
$$
Since the $2$'s in $a_i$ can only embed into $\getwo$'s in $w$, if $a_i$ is contained in a $k$-deletion of $w$ then its $1$ must embed into an element between $v(i-1)$ and $v(i)$, implying that $z(i)\ge 1$.  Conversely, if $a_i$ is not contained in a $k$-deletion of $w$ then either $\|a_i\|>n-k$ or $z(i)=0$.  Simple accounting shows that
$$
n-k=\left( (\mbox{\# of $1$'s in $w$}) + 2\ell + \ex_2(w)\right)-k,
$$
so $\|a_i\|=2\ell+1\le n-k$ because $w$ has more than $k$ $1$'s, and thus
\begin{equation}\label{eqn-det-z=0}
z(i)=0\iff\mbox{$a_i$ is not contained in a $k$-deletion of $w$.}
\end{equation}

The second step is to use these $0$'s to divine the nonzero entries of $z$.  Define the composition $b_i=b_i(1)\cdots b_i(\ell)$ by
$$
b_i(j)=
\left\{
\begin{array}{rll}
1&\mbox{if}&\mbox{$j\le i-1$ and $z(j)=0$ or}\\
&&\mbox{$j\ge i$ and $z(j+1)=0$, or}\\
2&\multicolumn{2}{l}{\mbox{otherwise,}}
\end{array}
\right.
$$
and consider the possible embeddings of $b_i$ in $w$.  Suppose for the sake of example that $i\ge 4$.  If $z(1)\ge 1$ then $b_i(1)=2$ and thus can embed only into or to the right of $v(1)$.  Otherwise if $z(1)=0$ then $b_i(1)=1$, but in this case $v(1)$ is the first entry of $w$ so again $b_i(1)$ can embed only into or to the right of $v(1)$.  Continuing this manner, if $z(2)\ge 1$ then $b_i(2)=2$, and since $b_i(2)$ can only embed into a $\getwo$ entry in $w$ to the right of $b_i(1)$, $b_i(2)$ can only embed into or to the right of $v(2)$.  Otherwise if $z(2)=0$ then $b_i(2)=1$, but then $v(1)$ and $v(2)$ are adjacent in $w$ so since $b_i(1)$ must embed into or to the right of $v(1)$ and $b_i(2)$ must embed to the right of $b_i(1)$ we see that $b_i(2)$ must embed into or to the right of $v(2)$.  Continuing in this manner it is easy to see (or more formally, to prove inductively) that:
\begin{itemize}
\item For all $j\le i-1$, $b_i(j)$ must embed into or to the right of $v(j)$.
\item For all $j\ge i$, $b_i(j)$ must embed into or to the left of $v(j)$.
\end{itemize}
These two facts combine to show that $b_i(i-1)$ and $b_i(i)$ can only embed between $v(i-1)$ and $v(i)$ (inclusive).  Now define the word $x\in\mathbb{N}^{\ell+1}$ by $x(i)=0$ if $z(i)=0$ and otherwise
$$
x(i)=\max\{s : \mbox{$b_i(1)\cdots b_i(i-1)\underbrace{1\cdots1}_{s}b_i(i)\cdots b_i(\ell)$ is contained in a $k$-deletion of $w$}\}.
$$
The analogue to \eqref{eqn-u-v-w} now follows by the conditions on embeddings of $b_i$ established above:
\begin{equation}\label{eqn-many-ones-zi}
x(i)=\min\{ z(i), n-k-\|b_i\| \}.
\end{equation}
Suppose $z(i)\ge 1$.  In this case $\|b_i\|=2\ell-h$, where $h$ denotes the number of $0$ entries of $z$ (``holes'').  Letting $k+t$ denote the number of $1$'s in $w$, we have
$$
n=k+t+2\ell+\ex_2(w).
$$
allowing us to rewrite \eqref{eqn-many-ones-zi} as
\begin{equation}\label{eqn-many-ones-zi2}
x(i)=\min\{ z(i), h+t+\ex_2(w) \}.
\end{equation}
If $\|v\|+\|x\|=n$ then we must have $z=x$ and thus have successfully reconstructed $w$.  By \eqref{eqn-many-ones-zi2}, this will happen if $z$ has no entries greater than $h+t+\ex_2(w)$.  Suppose, for the sake of contradiction, that this does not occur, i.e., that $z$ contains an entry greater than $h+t+\ex_2(w)$.  Then each of the other $(\ell+1-h)-1$ nonzero entries of $z$ correspond to at least one $1$ in $w$, and thus we have
$$
k+t
=
\mbox{\# of $1$'s in $w$}
\ge
(h+t+\ex_2(w)+1)+(\ell-h)
=
t+\ell+\ex_2(w)+1.
$$
However, this implies that
$$
2k\ge t+2\ell+\ex_2(w),
$$
so
$$
3k\ge (k+t)+2\ell+\ex_2(w)=n,
$$
and this contradiction completes the proof of both the lemma and Theorem~\ref{thm-comp-reconstruct}.
\end{proof}

\begin{example}
Suppose the reconstruction algorithm is given the set of $3$-deletions
$$
\{
1222, 2212,
11122, 11212, 11221, 12112, 12211,
111112, 111121, 111211, 112111,
1111111\}.
$$
of an unknown composition $w$ of $n=10$.  This set contains $1^{10-3}=1111111$ and every $3$-deletion in the set contains a $1$, so Lemma~\ref{lem-determine-ones} shows that $w$ has more than $k=3$ $1$'s.  Thus we follow the proof of Lemma~\ref{lem-many-ones}.  Each of the compositions with the fewest $1$'s, e.g., $2122$, give the $\getwo$ entries of $w$ in their correct order, $v=222$, so
$$
w=\underbrace{1\cdots1}_{z(1)}2\underbrace{1\cdots1}_{z(2)}2\underbrace{1\cdots1}_{z(3)}2\underbrace{1\cdots1}_{z(4)}.
$$
We then find the $0$ entries of $z$:
$$
\begin{array}{l}
\mbox{$z(1)\neq 0$ because $a_1=1222$ is contained in a $3$-deletion of $w$,}\\
\mbox{$z(2)=0$ because $a_2=2122$ is not contained in a $3$-deletion of $w$,}\\
\mbox{$z(3)\neq 0$ because $a_3=2212$ is contained in a $3$-deletion of $w$,}\\
\mbox{$z(4)=0$ because $a_4=2221$ is not contained in a $3$-deletion of $w$.}
\end{array}
$$
Now we build the word $x\in\mathbb{N}^4$.  We have that $x(2)=x(4)=0$ because the corresponding entries of $z$ are $0$.  To compute the other entries of $x$ we construct $b_1=121$ and $b_3=211$ and then have
$$
\begin{array}{l}
\mbox{$x(1)=3$ because $111\ 121$ is contained in a $3$-deletion of $w$ but $1111\ 121$ is not,}\\
\mbox{$x(3)=1$ because $21\ 1\ 1$ is contained in a $3$-deletion of $w$ but $21\ 11\ 1$ is not.}
\end{array}
$$
Since $\|v\|+\|x\|=\|222\|+\|3010\|=10$, we must have $z=x$ and thus $w=1112212$.
\end{example}

\minisec{The connection to permutations}
The subject of permutation patterns (see B\'ona's text~\cite{bona:combinatorics-o:} for a survey) is concerned with the following partial order on permutation: for permutations $\sigma$ of length $k$ and $\pi$ of length $n$, let $\sigma\le\pi$ if there are indices $i_1<i_2<\cdots<i_k$ such that the subsequence $\pi(i_1)\pi(i_2)\cdots\pi(i_k)$ has the same pairwise comparisons as $\sigma(1)\sigma(2)\cdots\sigma(k)$, and in such a case $\sigma$ is said to be an $(n-k)$-deletion of $\pi$.  For example, $13254\le 213654798$ because of the subsequence $26598$ ($=\pi(1)\pi(4)\pi(5)\pi(8)\pi(9)$).

Given two permutations $\sigma$ and $\pi$ of lengths $m$ and $n$ respectively, their {\it direct sum\/}, $\sigma\oplus\pi$, is the permutation of length $m+n$ whose first $m$ entries form $\sigma$ and whose last $n$ entries are the copy of $\pi$ obtained by adding $m$ to each entry.  For example, $213654\oplus 132=213654798$.  A permutation is said to be {\it layered\/} if it can be written as the direct sum of decreasing permutations.  Thus $213654798$ is layered because it can be written as $21\oplus 1\oplus 321\oplus 1\oplus 21$.  There is a natural order-preserving bijection between layered permutations and compositions; for example, $213654798=21\oplus 1\oplus 321\oplus 1\oplus 21$ maps to the composition $21312$ while $13254=1\oplus 21\oplus 21$ maps to $122$, and $122\le 21312$ under the partial order on compositions.

Smith~\cite{smith:permutation-rec:} was the first to study multiset reconstruction for permutations.  Her work was followed by Raykova~\cite{raykova:permutation-rec:} who proved that for all $k$, all sufficiently long permutations are reconstructible from their multisets of $k$-deletions.  This leaves open the question of whether all sufficiently long permutations are reconstructible from their {\it sets\/} of $k$-deletions.  Our work therefore answers Raykova's question about whether all sufficiently long layered permutations can be reconstructed from their sets of $k$-deletions.

\minisec{Acknowledgement}
I thank Robert Brignall for his helpful comments.

\bibliographystyle{acm}
\bibliography{../refs}

\end{document}